\magnification=1200
\overfullrule=0mm
\nopagenumbers

\centerline {\bf The Multivariate Fundamental Theorem of Algebra 
and Algebraic Geometry}
\medskip
\centerline {H. Hakopian}
\medskip

We derive two consequences of the multivariate fundamental theorem of algebra (MFTA)presented in [1-4]. The first one is the Bezout theorem for $n$ polynomials $g_1,\ldots,g_n$ from $k[x_1,\ldots,x_n].$ Notably the intersection multiplicities: $I_x({\cal G}),\ ({\cal G}=\{g_1,\ldots,g_n\})$ as in MFTA, are characterized just by means of partial differential operators given by polynomials from $D$-invariant linear spaces $D_x({\cal G})$ (see [3]):
$$D_x({\cal G}):=\{ p: \left[D^\alpha p\right](D)g_i(x)=0,\quad i=1,\ldots,n,\ \hbox{for all}\
\alpha \in Z_+^n\},$$
$$I_x({\cal G}):=\dim D_x({\cal G}).$$
\proclaim Theorem (Bezout). Let $k$ be an algebraically closed field.
Suppose that there is no intersection at ``infinity". 
Then 
$$ \sum_xI_x({\cal G})=
\cases
{\deg (g_1)\cdots \deg (g_n),&or\cr
	\infty.\cr}
$$
\par
\noindent The case of the intersection at ``infinity" can be treated in a standard way.

\noindent The second consequence is the following

\proclaim Theorem. Let $k$ be an algebraically closed field and
$f,f_1,\ldots,f_s \in k[x_1,\ldots,x_n].$ Then $f$ belongs to the polynomial ideal $<f_1,\ldots,f_s>$ if and only if for any $p\in k[x_1,\ldots,x_n],$
$$\left[D^\alpha p\right](D)f_i(x)=0,\quad i=1,\ldots,s,\ \hbox{for all}\
\alpha \in Z_+^n, \quad \hbox{implies}$$
$$\left[D^\alpha p\right](D)f(x)=0,\quad \hbox{for all}\
\alpha \in Z_+^n.$$
\par
\noindent Let us mention that one readily gets Nullstellensatz from here.
Moreover, the integer $m$ there, i.e., the power of $f,$ can be chosen such that
$m\le \deg(f_1)\cdots \deg (f_s)+1.$

\medskip
\noindent {\bf References}
\medskip
 
\item{1.} {\rm Hakopian, H.; Tonoyan, M.},
Polynomial interpolation and a multivariate analog of 
fundamental theorem of algebra, {\sl Symposium on Trends in Approximation Theory, Nashville} (2000), Abstracts, p.50, U.S.A.
\item{2.} {\rm Hakopian, H.; Tonoyan, M.},
Polynomial interpolation and a multivariate analog of the fundamental theorem of algebra, {\sl East J. on Approx.} {\bf 8} (2002) 355-379.
\item{3.} {\rm Hakopian, H.}, 
A multivariate analog of fundamental theorem of algebra 
and Hermite interpolation, {\sl Constructive Theory of Functions}, Ed.  Bojanov, B., Darba, Sofia, 2003, 1-18.
\item{4.} {\rm Mourrain, B.},
A new criterion for normal form algorithms, {\it in} {\sl Applied Algebra, Algebraic Algorithms and Error-Correcting Codes,} 13th Intern.\ Symp., AAECC-13, Honolulu, Hawaii, U.S.A., Nov.'99, Proc., Fossorier, M.;
Imai, H.; Lin, S.; Pol, A. (Eds.), Springer Lecture Notes in Computer Science, 1719, Springer-Verlag (Heidelberg), (1999), 430-443.

\bye